\documentclass[letterpaper, 10pt, conference]{ieeeconf}  

\IEEEoverridecommandlockouts                              
\usepackage{cite}
\usepackage{xargs}                      
\usepackage[noend]{algorithmic}
\usepackage{mathematix}
\usepackage{hyperref}
\usepackage{nicefrac}
\usepackage{microtype}
\usepackage{units}
\usepackage{flushend}

\hypersetup{
    bookmarks=true,         
    unicode=false,          
    pdftoolbar=true,        
    pdfmenubar=true,        
    pdffitwindow=false,     
    pdfstartview={FitH},    
    pdftitle={},    
    pdfauthor={},     
    pdfkeywords={},
    pdfsubject={Confidential},   
    pdfnewwindow=true,      
    colorlinks=true,        
    linkcolor=red,          
    citecolor=blue,         
    filecolor=magenta,      
    urlcolor=cyan           
}

\newcommand{\CLip}{C^{1{},{}1}}
\newcommand{\CLipL}[1]{\CLip_{{#1}}}

\newcommand{\ttop}{{\scriptstyle\top}}
\newcommand{\tttop}{{\scriptscriptstyle\top}}

\title{\LARGE \bf
Embedded nonlinear model predictive control for obstacle avoidance using PANOC
}

\author{Ajay~Sathya, Pantelis~Sopasakis, Ruben~Van~Parys, Andreas~Themelis, Goele~Pipeleers and Panagiotis~Patrinos
\thanks{$^{1}$A. Sathya, R. Van Parys and G. Pipeleers are with the Dept. of Mechanical Engineering,  \textsc{meco} research group, KU Leuven, 3001 Leuven, Belgium.}%
\thanks{$^{2}$P. Sopasakis, A. Themelis and P. Patrinos are with the Dept. of Electrical Engineering (\textsc{esat-stadius}), KU Leuven, 3001 Leuven, Belgium.}%
\thanks{The work of the these  authors was supported by the Research Foundation Flanders (FWO) PhD fellowship 1196818N; FWO research projects G086518N and G086318N; KU Leuven internal funding StG/15/043; Fonds de la Recherche Scientifique -- FNRS and the Fonds Wetenschappelijk Onderzoek -- Vlaanderen under EOS Project no 30468160 (SeLMA). R. Van Parys is a PhD fellow of the Research Foundation Flanders (FWO-Flanders).}%
        }

\begin{document}
\maketitle
\thispagestyle{empty}
\pagestyle{empty}
%
%
\begin{abstract}
We employ the proximal averaged Newton-type method for optimal control (PANOC) to solve obstacle avoidance problems in real time. We introduce a novel modeling framework for obstacle avoidance which allows us to easily account for generic, possibly nonconvex, obstacles involving polytopes, ellipsoids, semialgebraic sets and generic sets described by a set of nonlinear inequalities. PANOC is particularly well-suited for embedded applications as it involves simple steps, its implementation comes with a low memory footprint and its fast convergence meets the tight runtime requirements of fast dynamical systems one encounters in modern mechatronics and robotics. The proposed obstacle avoidance scheme is tested on a lab-scale autonomous vehicle.
\end{abstract}
%
%
\section{Introduction}\label{sec:intro}
\subsection{Background and Contributions}
Autonomous navigation in obstructed environments is a key element in emerging applications such as driverless cars, fleets of automated vehicles in warehouses and aerial robots performing search-and-rescue expeditions. Well-known approaches for finding collision-free motion trajectories are graph-search methods~\cite{takahashi1989voronoi}, virtual potential field methods~\cite{minguez2016motion_planning} or methods using the concept of velocity obstacles \cite{bareiss2015generalized}.

Recent motion planning research focuses on optimization-based strategies. Here, an optimal motion trajectory is sought while collision-avoidance requirements are imposed as constraints thereon. There is a broad range of such formulations. The most straightforward approach constrains the Euclidean distance between the vehicle and obstacle centers~\cite{wang2014synthesis}. This, however, only allows to separate spheres and ellipsoids as the computation of distances from arbitrary sets is generally not an easy operation.

Other approaches demand the existence of a hyperplane between the vehicle and obstacle at each time instant~\cite{mercy2017splines,debrouwere2013time}. This allows the separation of general convex sets.
Some methods formulate the collision avoidance requirement by mixed-integer constraints~\cite{alrifaee2014collision}. These types of problems are, however, cumbersome to solve in real-time.
For some specific problem types it is possible to transform the system dynamics such that collision-avoidance translates into simple box constraints on the states~\cite{frasch2013nmpc,turri2013lane}.

Model predictive control (MPC) is a powerful control strategy where control actions are computed by optimizing a cost function which is chosen in order to achieve a control task. Constraints on states, inputs and outputs can be seamlessly incorporated into such a framework. When the system dynamics is linear, the constraints are affine and the cost functions are quadratic, the associated optimization problem is a quadratic program. There exists a mature machinery of convex optimization algorithms~\cite{nocedal2006numerical,boyd2004convex} which are fast, robust and possess global convergence guarantees which can be used to solve these problems.

Nevertheless, the dynamics of most systems of interest are better modeled by nonlinear equations and constraints are often nonconvex. This situation is very common in obstacle avoidance involving nonconvex optimization problems. These are commonly solved using \textit{sequential quadratic programming} (SQP)~\cite{nocedal2006numerical} and \textit{interior point} (IP) methods~\cite{wachter2006intpoint} which are not well-suited for embedded applications with tight runtime requirements. Nonlinear MPC is often performed using the real-time iteration scheme proposed in~\cite{diehl2005realtime} which trades speed for accuracy and is accompanied by global convergence guarantees under certain assumptions.

For an algorithm to be suitable for an embedded implementation, it needs to involve only simple steps. This deems methods of the forward-backward-splitting (FBS) type~\cite{attouch2013fbs}, such as the proximal gradient method~\cite[Sec.~2.3]{bertsekas1999nonlinear}, appealing candidates. FBS-type algorithms can be used to solve nonlinear optimal control problems with simple input constraints via first eliminating the state sequence and expressing the cost as a function of the sequence of inputs alone --- the so-called \textit{single shooting} formulation. However, despite its simplicity, FBS, like all first-order methods, can exhibit slow convergence. Its convergence rate is at best Q-linear with a Q-factor close to one for ill-conditioned problems such as most nonlinear MPC problems.

In this work we propose a new modeling framework for generic constraints which can accommodate general nonconvex sets. The proposed methodology assumes that the obstacles are described by a set of nonlinear inequalities and does not require the computation of projections or distances to them. Then, the obstacle avoidance constraints are written as a nonlinear equality constraint involving a smooth function which, in turn, is relaxed using a penalty function. 

The resulting problems are solved using PANOC, a proximal averaged Newton-type method for optimal control, which was recently proposed in~\cite{stella2017simple}. Gradients of the cost function can be efficiently computed using automatic differentiation toolboxes such as CasADi. The algorithm is simple to implement, yet robust, since it combines projected gradient iterations with quasi-Newtonian directions to achieve fast convergence.

The proposed framework is tested on a number of simulation scenarios where we show that it is possible to avoid obstacles of complex shape described by nonlinear inequalities. PANOC is compared with SQP and IP methods and is found to be significantly faster. Furthermore, we present experimental results on a lab-scale robotic platform which runs a C implementation of PANOC.

\subsection{Notation}
Let $\N$ be the set of nonnegative integers, $\N_{[k_1, k_2]}$ be the set of integers in the interval $[k_1, k_2]$ and $\barre = \Re\cup\{+\infty\}$ be the set of extended real numbers. For a matrix $A\in\Re^{m\times n}$, we denote its transpose by $A^\ttop$. For $x\in\Re$, we define the operator $\plus{x} = \max\{x,0\}$.
For a nonempty closed convex set $U\subseteq \Re^n$, the \textit{projection} onto $U$ is the operator $\proj_U(v) = \argmin_{u\in U} \|u-v\|$. The \textit{distance} from the set $U$ is defined as $\dist_U(v) = \inf_{u\in U}\|u-v\|$.
The class of continuously differentiable functions $f{}:{}\Re^n\to\Re$ is denoted by $C^1$. The subset of $C^1$ of functions with Lipschitz-continuous gradient is denoted as $\CLip$. We use the notation $\CLipL{L}$ for $\CLip$ functions with $L$-Lipschitz gradients.

\section{NMPC for obstacle avoidance}
\label{sec:obstacle_avoidance}

\subsection{Problem statement}

Kinematic equations lead to continuous-time nonlinear dynamical systems of the form $\dot{x} = f_c(x, u, t)$ where $x\in\Re^{n_x}$ is the system state, typically a vector comprising of position, velocity and orientation data, and $u\in\Re^{n_u}$ is the control signal. We assume that the position coordinates $z\in\Re^{n_d}$ are part of the state vector. The continuous-time dynamics can be discretized (for instance, using an explicit Runge-Kutta method) leading to a discrete-time dynamical system of the form
\begin{equation}
 x_{k+1} = f_k(x_k, u_k).
\end{equation}
As it is typically the case in practice, we assume that $f_k:\Re^{n_x}\times \Re^{n_u}\to \Re^{n_x}$ are smooth mappings.

The objective of the navigation controller is to steer the controlled vehicle from an initial state $x_0$ to a target state $x_{\mathrm{ref}}$, typically a position in space together with a desired orientation. At the same time, the vehicle has to avoid certain, possibly moving, obstacles which are described by open sets $O_{kj} \subseteq \Re^{n_d}$, $j\in\N_{[1,q_k]}$,
each described by
\begin{equation}
 O_{kj} = \{z\in\Re^{n_d} {}:{} h_{kj}^i(z) > 0, i\in\N_{[1,m_{kj}]}\}.
\end{equation}
Sets $O_{kj}$ need not be convex.
Obstacle avoidance constraints can be concisely written as 
\begin{equation}
\label{eq:oa_constraints}
 z_k\notin O_{kj},\text{ for } j\in\N_{[1,q_k]}.
\end{equation}
Moreover, the vehicle is only allowed to move in a domain which is described by the inequalities
\begin{equation}
 \label{eq:domain_constraints}
 g_k(x_k, u_k) \leq 0,
\end{equation}
where $g_k:\Re^{n_x}\times \Re^{n_u}\to \Re^{n_c}$ is a $C^2$ mapping and $\leq$ is meant in the element-wise sense.
Control actions $u_k$ are constrained in a closed compact set $U_k$, that is
\begin{equation}
\label{eq:Uk_set}
 u_k \in U_k,
\end{equation}
on which it is easy to project and hereafter shall be assumed to be convex. Sets $U_k$ often represent box constraints of the form $U_k = \{u\in\Re^{n_u} {}:{} u_{\mathrm{min}} \leq u \leq u_{\mathrm{max}}\}$.

\subsection{Nonlinear model predictive control}
\label{subsec: nmpc cost}
Nonlinear model predictive control problems arising in obstacle avoidance can be written in the following form
\begin{subequations}\label{eq:nmpc-1}
\begin{align}
 \minimize\quad &\ell_N(x_N) + \sum_{k=0}^{N-1}\ell_k(x_k, u_k), \label{eq:nmpc-1-1}\\
 \stt\quad &x_{0} = x,\label{eq:nmpc-1-2}\\
 &x_{k+1} = f_k(x_k, u_k),\, k\in\N_{[0,N-1]},\label{eq:nmpc-1-3}\\
 &u_k \in U_k,\, k\in\N_{[0,N-1]},\label{eq:nmpc-1-4}\\
 &z_k\notin O_{kj},\, j\in\N_{[1,q_k]},\, k\in\N_{[0,N]}\label{eq:nmpc-1-5}\\
 &g_k(x_k, u_k) \leq 0,\, k\in\N_{[0,N]},\label{eq:nmpc-1-6}\\
 &g_N(x_N) \leq 0.\label{eq:nmpc-1-7}
\end{align}
\end{subequations}
The stage costs $\ell_k:\Re^{n_x}\times \Re^{n_u}\to \Re$ for $k\in\N_{[0,N-1]}$ in~\eqref{eq:nmpc-1-1} are $\CLip$ functions penalizing deviations of the state from the reference (destination and orientation) and may be taken to be quadratic functions of the form 
$\ell(x_k, u_k) = 
(x_k-x_{\mathrm{ref}})^{\tttop} Q_k (x_k-x_{\mathrm{ref}}) + 
(u_k-u_{\mathrm{ref}})^{\tttop} R_k (u_k-u_{\mathrm{ref}})$. The terminal cost $\ell_N:\Re^{n_x}\to \Re$ in~\eqref{eq:nmpc-1-1} is a $\CLip$ function such as $\ell_N(x_N) = (x_N - x_{\mathrm{ref}})^{\tttop} Q_N (x_N - x_{\mathrm{ref}})$.

\subsection{Reformulation of obstacle avoidance constraints}
\label{sec:reform-obst-avoid}
Consider an obstacle described as the intersection of a finite number of strict nonlinear inequalities 
\begin{equation}
\label{eq:obstacle}
O = \{z\in\Re^{n_d} {}:{} h^i(z) > 0, i\in\N_{[1,m]}\}, 
\end{equation}
where $h^i:\Re^{n_d} \to \Re$ are $\CLip$ functions.
 The constraint $z\notin O$ --- cf.~\eqref{eq:nmpc-1-5} --- is satisfied if and only if 
 \begin{equation}
  h^{i_0}(z) \leq 0,\,\text{for some}\, i_0\in\N_{[1,m]},
 \end{equation}
or, equivalently, $\plus{h^{i_0}(z)}^2=0$. This constraint can then be encoded as
\begin{equation}
\psi_O(z) \dfn \tfrac{1}{2}\prod_{i=1}^{m} \plus{h^{i}(z)}^2 = 0.\label{eq:obst-avoidance-equality}
\end{equation}
We have expressed the obstacle avoidance constraints as a nonlinear equality constraint.
We should remark that, unlike approaches based on the distance-to-set function~\cite{gilbert1985distance}, function $\psi_O$ in~\eqref{eq:obst-avoidance-equality} is a $C^1$ function of $z$. Indeed, $\psi_O$ is differentiable in $\Re^{n_d}\setminus O$ (and equal to $0$), it is differentiable in $O$ with gradient
\begin{equation}
 \label{eq:psi-gradient}
 \nabla \psi_O(z) = 
 \begin{cases}
  \sum\limits_{i{}={}1}^{m} \hspace*{-1pt}h^i(z)\hspace*{-2pt}\prod\limits_{j{}\neq{}i}\! (h^j(z))^2 \nabla h^i(z),&\hspace*{-5pt}\text{if}\,x{}\in{}O,\\
  0,&\hspace*{-5pt}\text{otherwise.}
 \end{cases}
\end{equation}
Note that $\nabla \psi_O$ on the boundary of $O$ vanishes, so it is everywhere continuous.
If, additionally, $O$ is bounded, $\psi_O$ is $C^{1,1}$.

The formulation of Eq.~\eqref{eq:obst-avoidance-equality} can be used for obstacles described by quadratic constraints of the form 
$
O = \{z\in\Re^{n_d}{}:{} 1- (z-c)^\tttop E (z-c) > 0\}, 
$
such as balls and ellipsoids. Then, the associated equality constraint becomes
$
 \plus{1 - (z-c)^\tttop E (z-c)}^2 = 0.
$
Polyhedral obstacles of the form 
$
O=\{z\in\Re^{n_d}{}:{} b_i-a_i^{\tttop}z > 0,\in\N_{[1,m]}\}, 
$
with $b_i\in\Re$ and $a_i\in\Re^{n_d}$, can also be accommodated by the constraint
$
 \prod_{i=1}^{m}\plus{b_i - a_i^\tttop z}^2 = 0.
$

Equation~\eqref{eq:obst-avoidance-equality} can also be used to describe general nonconvex constraints such as ones described by semi-algebraic sets where $h^i(z)$ are polynomials as well as any other obstacle which is available in the aforementioned representation.

Obstacle avoidance constraints~\eqref{eq:obst-avoidance-equality} are equivalent to the existence of $t\in\Re^m$ so that
\begin{subequations}
\begin{align}
 \min \{t_1,\ldots, t_m\} = 0,\\
 h^i(z) \leq t_i, \text{ for } i\in\N_{[1,m]}.
\end{align}
\end{subequations}
This observation reveals a link to vertical complementarity constraints which have been studied extensively in the literature~\cite{scheel2000complem,liang2012vcc}.

\subsection{Relaxation of constraints}
\label{sec:relaxation}
Due to the fact that modeling errors and disturbances may lead to the violation of imposed constraints and infeasibility of the MPC optimization problem, it is common practice in MPC to replace state constraints by appropriate penalty functions known as \textit{soft constraints}. Quadratic penalty functions are often used for this purposes revealing a clear link between this approach and the \textit{quadratic penalty method} in numerical optimization~\cite[Sec.~4.2.1]{bertsekas1999nonlinear}.

Equality constraints, such as the ones arising in the reformulation of the obstacle avoidance constraints in Section~\ref{sec:reform-obst-avoid}, can be relaxed by means of soft constraints. Indeed, constraints of the form $\Phi_{k}(z_k) = 0$, $k\in\N_{[1,N]}$, where $\Phi_k:\Re^{n_d}\to\Re_+$ are $C^2$ functions, can be relaxed by introducing the penalty functions $\tilde{\Phi}_k(z_k) = \eta_k \Phi_{k}(z_k)$ for some weight factors $\eta_k>0$. 

That said, constraints like~\eqref{eq:obst-avoidance-equality} for a set of time-varying obstacles $O_{kj}=\{z\in\Re^{n_d} {}:{} h_{kj}^i(z) > 0, i\in\N_{[1,m_{kj}]}\}$, with $j\in\N_{[1,q_k]}$, can be relaxed by the appending the following term in the original cost function
\begin{equation}
 \label{eq:soft_penalty}
 \tilde{h}_{k}(z_k) = \sum_{j=1}^{q_k}\eta_{kj}\prod_{i=1}^{m_{kj}} \plus{h^{i}_{kj}(z)}^2,
\end{equation}
for some positive weight factors $\eta_{kj}>0$.

Note that the proposed approach for dealing with obstacle avoidance constraints requires only a representation of the obstacles in the generic form~\eqref{eq:obstacle} and does not call for the computation of distances to the obstacles as in distance-based methods~\cite{wang2014synthesis,gilbert1985distance}, nor does it require the obstacles to be convex sets.

Similarly, inequality constraints of the form $g_k(x_k, u_k)\leq 0$ and $g_N(x_N)\leq 0$ can be relaxed by introducing the penalty functions $\tilde{g}_k(x,u) = \beta_{k}\plus{g_k(x, u)}^2$ and $\tilde{g}_N(x) = \beta_N \plus{g_N(x)}^2$ for positive weights $\beta_k>0$, $k\in\N_{[0,N]}$.

We may now relax the state constraints in~\eqref{eq:nmpc-1} by defining the modified stage cost and terminal cost functions
\begin{equation*}
\begin{aligned}
 \tilde{\ell}_k(x, u) &= \ell_k(x, u) + \tilde{g}_k(x, u)+ \tilde{h}_{k}(z_k),\\
 \tilde{\ell}_N(x) &= \ell_N(x) + \tilde{g}_N(x),
\end{aligned}
\end{equation*}
leading to the following relaxed optimization problem without state constraints
\begin{subequations}\label{eq:nmpc-2}
\begin{align}
 \minimize\quad &\tilde{\ell}_N(x_N) + \sum_{k=0}^{N-1}\tilde{\ell}_k(x_k, u_k), \label{eq:nmpc-2-1}\\
 \stt\quad &x_{0} = x,\label{eq:nmpc-2-2}\\
 &x_{k+1} = f_k(x_k, u_k),\, k\in\N_{[0,N-1]},\label{eq:nmpc-2-3}\\
 &u_k \in U_k,\, k\in\N_{[0,N-1]},\label{eq:nmpc-2-4} 
\end{align}
\end{subequations}

\subsection{NMPC problem formulation}
In this section we cast the nonlinear MPC problem~\eqref{eq:nmpc-1} as
\begin{equation}
 \minimize\limits_{u\in{} U}\, \ell(u), 
 \label{eq:nmpc-ss-1}
\end{equation}
where the optimization is carried out over vectors $u = (u_0, \ldots, u_{N-1})\in \Re^{n}$, with $n=Nu_u$ and $U\dfn U_0\times U_1\times\cdots\times U_{N-1}$ and $\ell:\Re^n\to\Re$ is a real-valued $\CLipL{L_\ell}$ function.

We introduce the following sequence of functions $F_k:\Re^n\to\Re^{n_x}$ for $k\in\N_{[0,N-1]}$
\begin{subequations}
\begin{align}
 F_0(u)     &= x,\\
 F_{k+1}(u) &= f_k(F_k(u), u_k).
\end{align}
\end{subequations}
Define the smooth function $\ell:\Re^n\to\Re$
\begin{equation}\label{eq:cost-ss}
\ell(u) \dfn \tilde{\ell}_N(F_N(u)) +
 \sum_{k=0}^{N-1}\tilde\ell_k(F_k(u), u_k) .
\end{equation}
The gradient of function $\ell$ in~\eqref{eq:cost-ss} can be computed by means of the reverse mode of automatic differentiation (also known as adjoint method or backpropagation) as shown in Alg.~\ref{alg:ad}~\cite{dunn1989efficient}.
\begin{algorithm}
\caption{Automatic differentiation for $\ell$ in~\eqref{eq:cost-ss}}\label{alg:ad}
\begin{algorithmic}[1]
 \REQUIRE $x_0\in\Re^{n_x}$, $u\in\Re^{n}$.%
 \ENSURE  $\ell(u)$, $\nabla \ell(u)$
 \STATE $\ell(u) {}\leftarrow{}0$
 \FOR{\(k=0,\ldots,N-1\)}
      \STATE $x_{k+1} {}\leftarrow{} f_k(x_k, u_k)$, $\ell(u) {}\leftarrow{} \ell(u) + \tilde\ell_k(x_k, u_k)$
 \ENDFOR{}
 \STATE $\ell(u) {}\leftarrow{} \ell(u) + \tilde\ell_N(x_N)$, $p_N {}\leftarrow{} \nabla \tilde\ell_N(x_N)$
 \FOR{\(k=N-1,\ldots,0\)}
      \STATE $p_k {}\leftarrow{} \nabla_{x_k}f_k(x_k, u_k)p_{k+1} + \nabla_{x_k}\tilde\ell_k(x_k, u_k)$
      \STATE $\nabla_{u_k} \ell_k(u) {}\leftarrow{} \nabla_{u_k}f_k(x_k, u_k)+\nabla_{u_k}\tilde\ell_k(x_k, u_k)$
 \ENDFOR{}
\end{algorithmic}
\end{algorithm}

%
%
%
Problem~\eqref{eq:nmpc-ss-1} is in a form that allows the application of the projected gradient iteration
\begin{equation}
 u^{\nu+1} = T_\gamma(u^\nu) \dfn \proj_{U}(u^\nu - \gamma \nabla \ell(u^\nu)),
 \label{eq:fbs_recursion}
\end{equation}
with $\gamma>0$. In particular, if $\ell\in\CLipL{L_\ell}$ and $\gamma < \nicefrac{2}{L_\ell}$, then all accumulation points of~\eqref{eq:fbs_recursion}, $u^\star$, are fixed points of $T_\gamma$ called \textit{$\gamma$-critical points}, that is~\cite[Prop.~2.3.2]{nocedal2006numerical}
\begin{equation}
 u^\star = T_\gamma(u^\star).
 \label{eq:fixed-point}
\end{equation}

\section{Fast nonlinear MPC}
The problem of finding a fixed point for $T_\gamma$ can be reduced to the equivalent problem of finding a zero of the \textit{fixed-point residual} operator which is defined as the operator
\begin{equation}
 R_\gamma(u) = \tfrac{1}{\gamma}(u - T_\gamma(u)).
\end{equation}
This motivates the adoption of a Newton-type iterative scheme 
\begin{equation}
 \label{eq:quasi_newton}
 u^{\nu+1} = u^{\nu} - H_{\nu} R_\gamma(u^\nu),
\end{equation}
where $H_\nu$ are invertible linear operators, appropriately chosen so as to encode first-order information about $R_\gamma$. This is done by enforcing the \emph{inverse secant condition} $s^\nu = H_{\nu+1} y^\nu$, for $s^\nu=u^{\nu+1}-u^{\nu}$ and $y^\nu = R_\gamma(u^{\nu+1})-R_\gamma(u^{\nu})$ and can be obtained by quasi-Newtonian methods such as the limited-memory BFGS (L-BFGS) method~\cite{nocedal2006numerical} which is free from matrix operations, requires only a limited number of inner products and is suitable for embedded implementations.
The main weakness of this approach is that convergence is only guaranteed in a neighborhood of a $\gamma$-critical point $u^\star$. We shall describe a \textit{globalization} procedure which hinges on the notion of the \textit{forward-backward envelope} function.

\subsection{The Forward-Backward envelope function}
\begin{figure}[t]
 \centering
 \includegraphics[width=0.75\columnwidth]{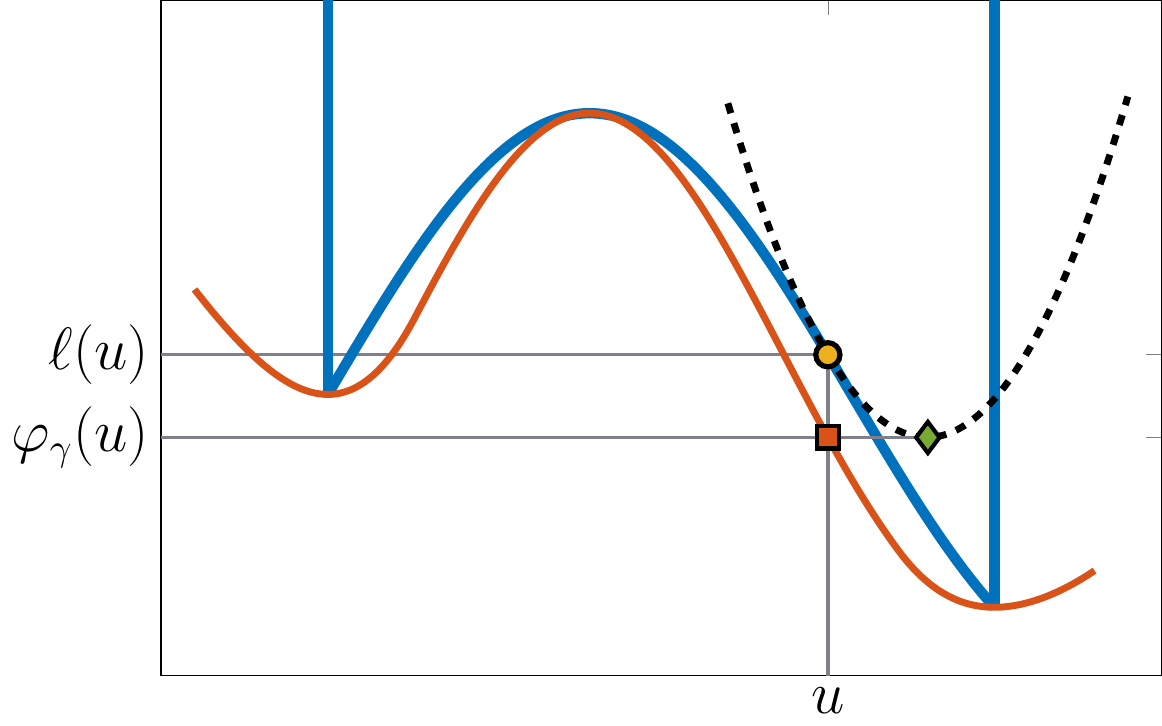}%
 \caption{%
      Construction of the FBE with $\gamma = 0.15$ (red line) for the function $\ell(u) = \sin(2u)$ over the set $U=[0,2]$ (thick blue line). The dashed line shows the approximation of $\ell$ at a point $u$ by the quadratic model $Q^\ell_\gamma(v;u)$. The value of FBE at $u$ is then given by $\varphi_\gamma(u) = \inf_{v\in U}Q^\ell_\gamma(v;u)$. Note that the local minima of $\ell$ over $U$ are exactly the local minima of $\varphi_\gamma$.%
      \label{fig:fbe}
      }%
\end{figure}

The \textit{forward-backward envelope} (FBE) is an exact, continuous, real-valued merit function for~\eqref{eq:nmpc-ss-1}~\cite{patrinos2014forward,stella2017forward,stella2017simple,themelis2016forward,PatBem13}. Function $\ell$ can be approximated at a point $u\in U$ by the quadratic upper bound
\begin{equation}
 Q^\ell_\gamma(v; u) = \ell(u) + \nabla \ell(u)^\tttop (v-u) + \nicefrac{1}{2\gamma}\|u-v\|^2.
\end{equation}
The FBE is then defined as 
\begin{equation}
 \varphi_\gamma(u) = \inf_{v\in U} Q^\ell_\gamma(v; u).
\end{equation}
This construction is illustrated in Fig.~\ref{fig:fbe}.
Provided that it is easy to compute the distance to $U$, the FBE can be easily computed by
\begin{equation}
 \varphi_\gamma(u) = \ell(u) - \tfrac{\gamma}{2}\|\nabla \ell(u)\|^2 + \tfrac{1}{2\gamma}\dist_U^2(u - \gamma \nabla \ell(u)).
\end{equation}
Therefore, the computation of the FBE is of the same complexity as that of a forward-backward step.

The FBE possesses several favorable properties, perhaps the most important being that it is real-valued, continuous and for $\gamma\in(0,\nicefrac{1}{L_\ell})$ shares the same (local/strong) minima with the original problem~\eqref{eq:nmpc-ss-1}.  This means that~\eqref{eq:nmpc-ss-1} is reduced to an unconstrained minimization problem. Moreover, if $\ell\in C^2$, then $\varphi_\gamma \in C^1$ with $\nabla \varphi_\gamma(u) = (I-\gamma \nabla^2 \ell(u))R_\gamma(u)$.

\subsection{PANOC Algorithm}
We employ the proximal averaged Newton-type method for optimal control, for short PANOC, which was recently proposed in~\cite{stella2017simple}. PANOC performs fast Newton-type updates while an FBE-based linesearch endows it with global convergence properties while it uses the same oracle as the projected gradient method.
\begin{algorithm}
\caption{PANOC algorithm for problem~\eqref{eq:nmpc-ss-1}}\label{alg:panoc}%
\begin{algorithmic}[1]
 \REQUIRE
 	\(\gamma\in(0,\nicefrac{1}{L_\ell})\),\
 	\(L_\ell>0\),\
 	\(\sigma\in(0,\tfrac{\gamma}{2}(1-\gamma L_\ell))\),\
 	\(u_0\in\Re^{n}\),\
 	\(x_0\in\Re^{n_x}\),
 	L-BFGS memory length $\mu$.%
\FOR{\(\nu=0,1,\ldots\)}
	\STATE\label{state:zerofpr2:FB}
		$\bar u^\nu{}\leftarrow{}\proj_{U}(u^\nu-\gamma\nabla\ell(u^\nu) )$
	\STATE	\label{state:zerofpr2:r}$r^\nu{}\leftarrow{}\gamma^{-1}(u^\nu-\bar u^\nu)$
 	\STATE \label{state:zerofpr2:d} \(d^\nu\leftarrow-H_\nu r^\nu\) using L-BFGS%
	\STATE \label{state:zerofpr2:upd} \(u^{\nu+1}{}\leftarrow{}u^\nu - (1-\tau\!_\nu)\gamma r^\nu + \tau\!_\nu d^\nu\),
		where \(\tau\!_\nu\) is the largest number in \(\{\nicefrac{1}{2^i}{}:{}i\in\N\}\) such that
		\begin{equation}\label{eq:LS}
			\varphi_\gamma(u^{\nu+1})
		{}\leq{}
			\varphi_\gamma(u^\nu)-\sigma\|r^\nu\|^2
		\end{equation}
\ENDFOR{}%
\end{algorithmic}
\end{algorithm}

The iterative scheme, which is presented in Alg.~\ref{alg:panoc}, involves the computation of a projected gradient point $\bar{u}^{\nu}$ in step~\ref{state:zerofpr2:FB} and an L-BFGS direction $d^{\nu}$ in step~\ref{state:zerofpr2:d}. L-BFGS obviates the need to store or explicitly update matrices $H_\nu$ in~\eqref{eq:quasi_newton} by storing a number $\mu$ of past values of $s^\nu$ and $y^\nu$. The computation of $d^\nu$ requires only inner products which amount to a maximum of $4\mu n$ scalar multiplications. In step~\ref{state:zerofpr2:upd}, the iterates are updated using a convex combination of the projected gradient update direction $-\gamma R_\gamma(u^\nu)$ and a fast quasi-Newtonian direction $d^{\nu}$. The algorithm is terminated when $\|R_\gamma(u^\nu)\|_{\infty}$ drops below a specified tolerance.

In line~\ref{state:zerofpr2:upd}, the backtracking line search procedure ensures that a sufficient decrease condition is satisfied using the FBE as a merit function. Under mild assumptions, eventually only fast updates are activated and updates reduce to  $u^{\nu{}+{}1} = u^{\nu} {}+{} d^{\nu}$. The sequence of fixed-point residuals $\{r^\nu\}_{\nu\in\N}$ converges to $0$ square summably, while PANOC produces sequences of iterates, $\{u^\nu\}_{\nu\in\N}$ and $\{\bar{u}^\nu\}_{\nu\in\N}$, whose cluster points are $\gamma$-critical points.

In absence of a Lipschitz constant $L_\ell$, the algorithm can be initialized with an estimate, e.g.,
$
 L_\ell^0 = {\|\nabla_u \ell(u^0 + \delta u) - \nabla_u \ell(u^0) \|}/{\|\delta u\|},
$
where $\delta u \in \Re^n$ is a small perturbation vector. Then, step~\ref{state:zerofpr2:FB} in Alg.~\ref{alg:panoc} needs to be replaced by the backtracking procedure of Alg.~\ref{alg:gamma_backtracking} which updates $L_\ell$, $\sigma$ and $\gamma$.
\begin{algorithm}[H]
\caption{Lipschitz constant backtracking}
\label{alg:gamma_backtracking}
\begin{algorithmic}
\WHILE{$\ell(\bar{u}^\nu) > \ell(u^\nu) - \gamma \nabla \ell(u^\nu)^\tttop r^\nu  + \nicefrac{L_\ell}{2}\|\gamma r^\nu\|^2$}
\STATE $L_\ell\leftarrow 2L_\ell$,\,
       $\sigma \leftarrow \nicefrac{\sigma}{2}$,\,
       $\gamma \leftarrow \nicefrac{\gamma}{2}$
\STATE $\bar u^\nu{}\leftarrow{}\proj_{U}(u^\nu-\gamma\nabla\ell(u^\nu) )$
\ENDWHILE
\end{algorithmic}
\end{algorithm}
Alg.~\ref{alg:gamma_backtracking} updates $L_\ell$, $\sigma$ and $\gamma$ only a finitely many times, so, it does not affect the convergence properties of the algorithm. 

In MPC, we may warm start the algorithm using the optimal solution of the previous MPC instance.
Unlike SQP and IP methods, PANOC does not require the solution of linear systems or quadratic programming problems at every iteration and it involves only very simple operations such as vector additions, scalar and inner products. Additionally, PANOC converges \textit{globally}, that is, from any initial point $u^0\in\Re^{n}$.

\section{Simulations}
\begin{figure}[t]
\centering
\includegraphics[width=0.7\linewidth]{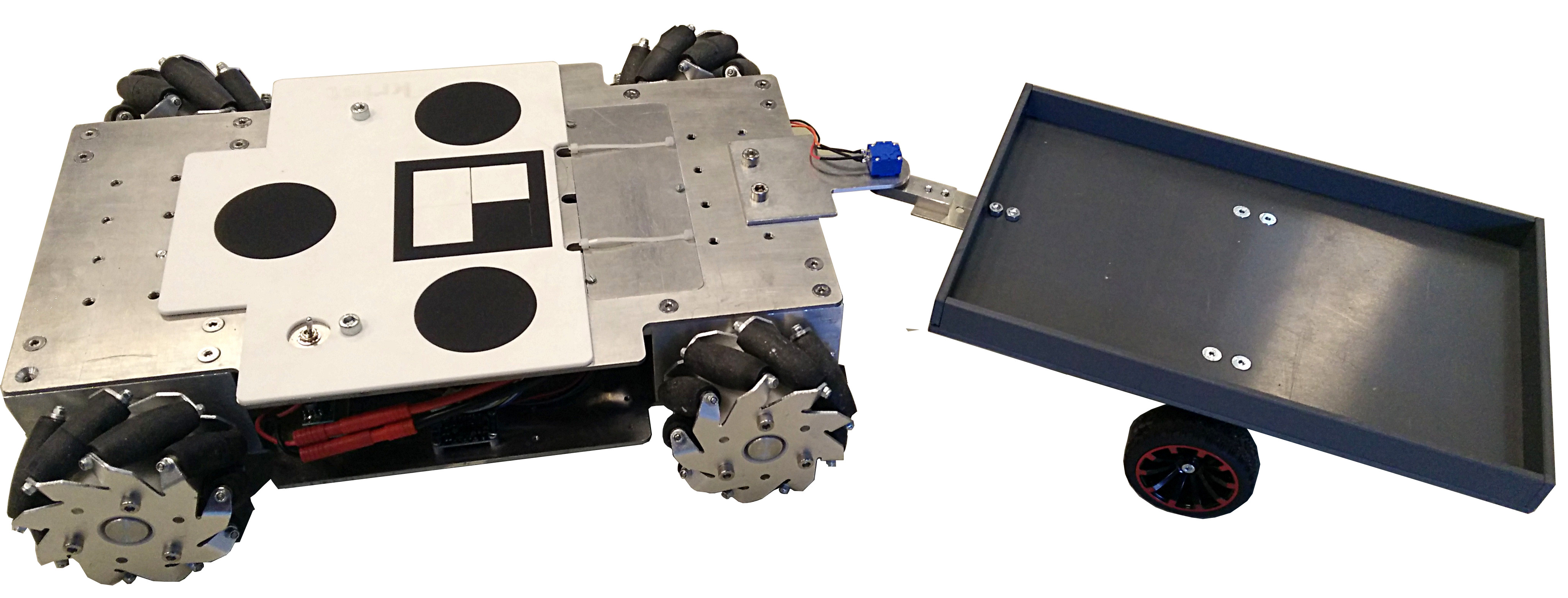}
\caption{In-house developed mobile robot with trailer.}
\label{fig:ourbot}
\end{figure}

\begin{figure}[t]
 \centering
 \includegraphics[width=0.47\columnwidth]{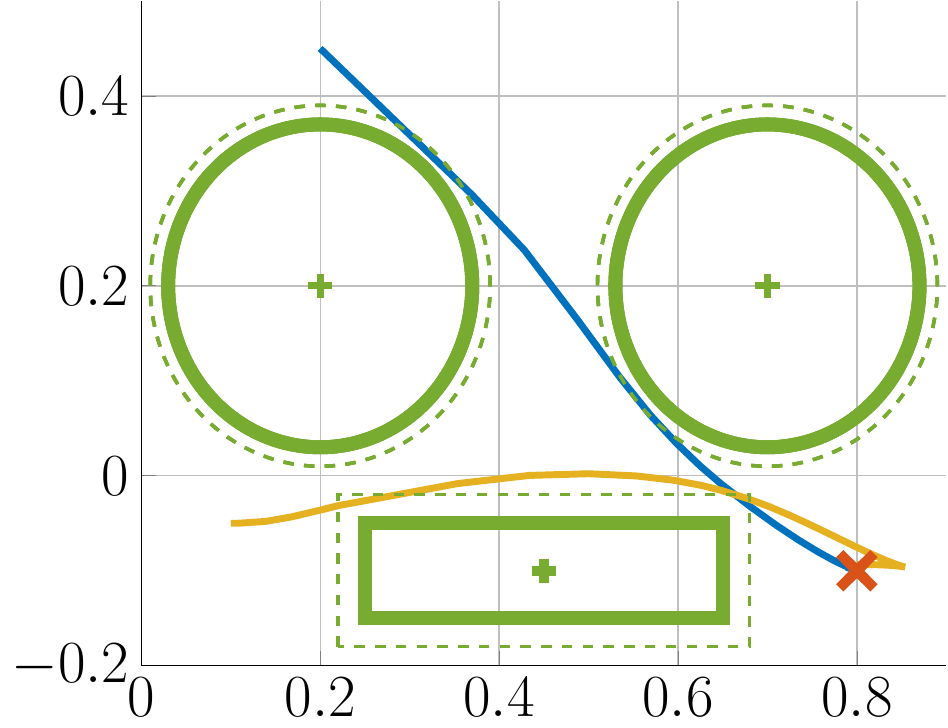}%
 \includegraphics[width=0.47\columnwidth]{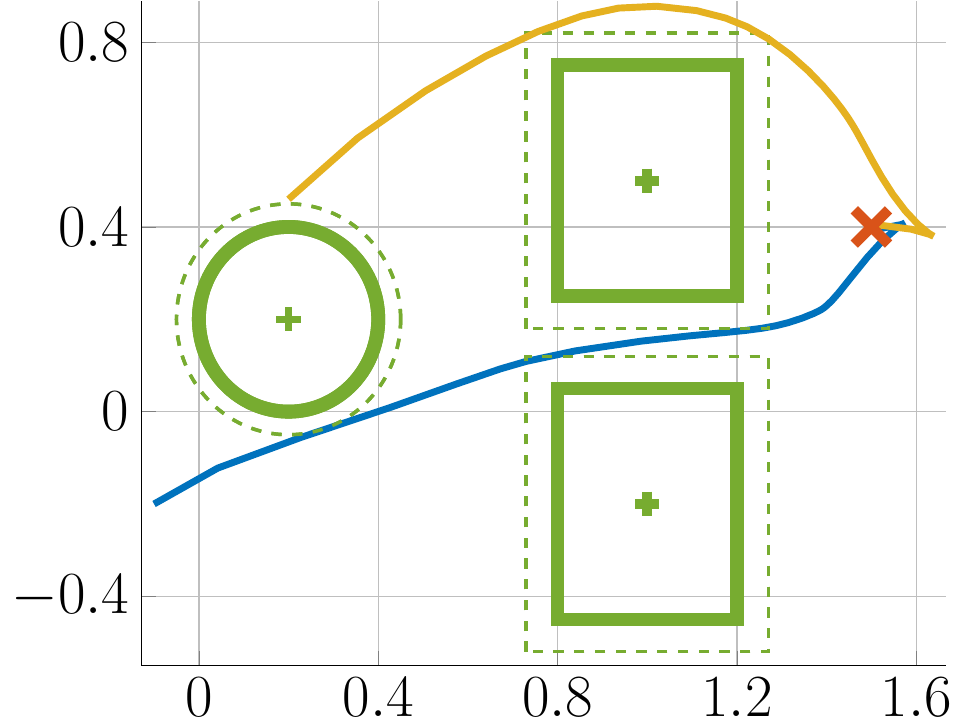}\\[1em]
 \includegraphics[width=0.47\columnwidth]{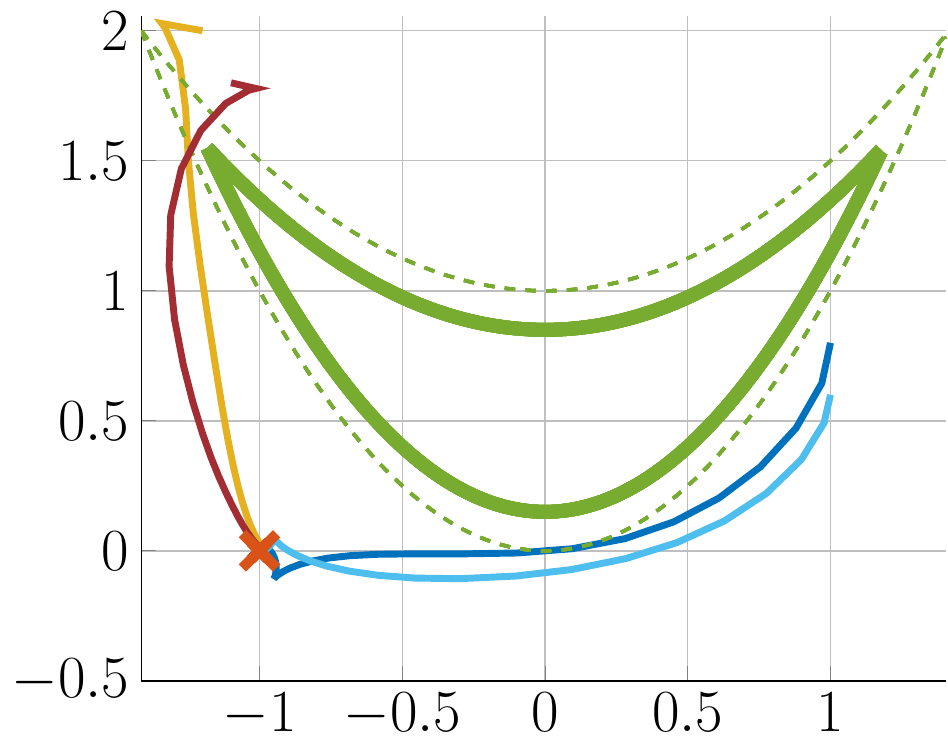}%
 \includegraphics[width=0.47\columnwidth]{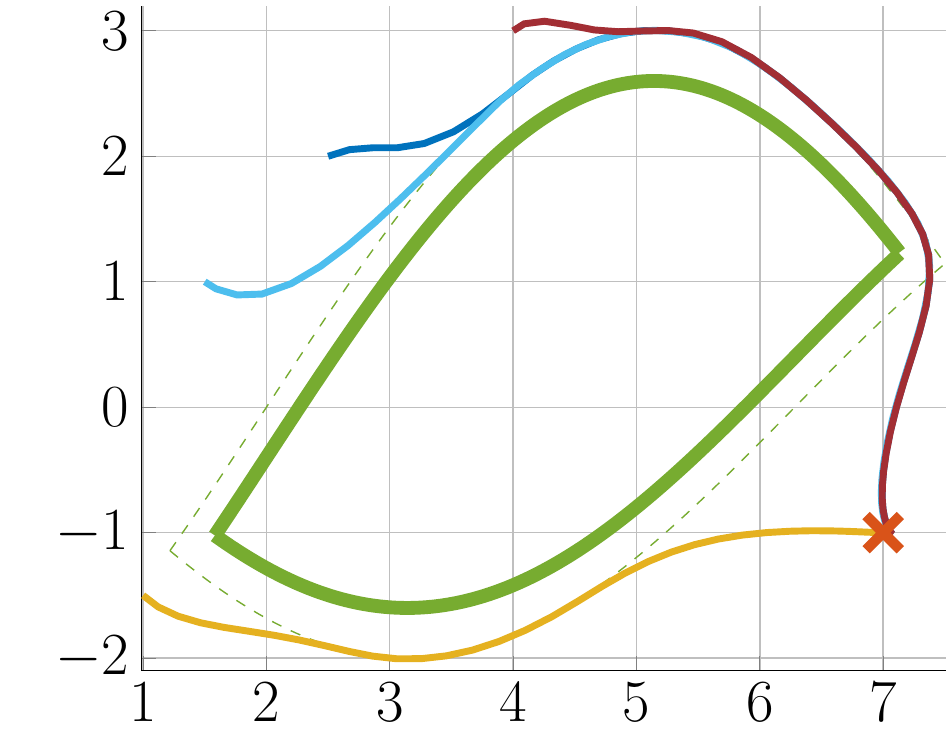}
 \caption{%
      Four obstacle avoidance scenarios using PANOC. The red \textsc{x} mark denotes the destination point and various lines are trajectories from different initial points. The dashed lines correspond to enlargements of obstacles which are  circumscribed by thick green lines. (First row) Obstacle avoidance with circular and rectangular obstacles, (Down, left) Enlarged obstacle defined by $O=\{(x,y){}:{} y > x^2, y < 1 + \nicefrac{x^2}{2}\}$ and (Down, right) Enlarged obstacle defined by $O = \{(x,y){}:{}  y > 2\sin(-\nicefrac{x}{2}), y < 3\sin(\nicefrac{x}{2}-1), 1 {<} x {<} 8\}$.
      \label{fig:obstacle_avoidance_scenarios}
      }%
\end{figure}

The proposed methodology is validated on the control of a mobile robot carrying a trailer shown in Fig.~\ref{fig:ourbot}. Assuming zero slip of the trailer wheels, the nonlinear kinematics is
\begin{subequations}\label{eq:kin_model}
\begin{align}
\dot{p}_x &= u_x + L\sin\theta {}\cdot{} \dot{\theta},\\
\dot{p}_y &= u_y - L\cos\theta{}\cdot{}\dot{\theta},\\
\dot{\theta} &= \tfrac{1}{L}\left(u_y\cos\theta - u_x\sin\theta \right),
\end{align}
\end{subequations}
where the state vector $x=(p_x, p_y, \theta)$ comprises the coordinates $p_x$ and $p_y$ of the trailer and the heading angle $\theta$. The input $u=(u_x, u_y)$ is a velocity reference which is tracked by a low-level controller. The distance between the center of mass of the trailer and the fulcrum connecting to the towing vehicle is $L = \unit[0.5]{m}$.
In Fig.~\ref{fig:obstacle_avoidance_scenarios} we present four obstacle avoidance scenarios involving, among other, obstacles described by polynomial and trigonometric functions. The system dynamics given in~\eqref{eq:kin_model} is discretized using the fourth-order Runge-Kutta method and we use memory $\mu=10$ for the L-BFGS directions in Alg.~\ref{alg:panoc}.

\begin{figure}[t]
 \centering
 \includegraphics[width=0.8\columnwidth]{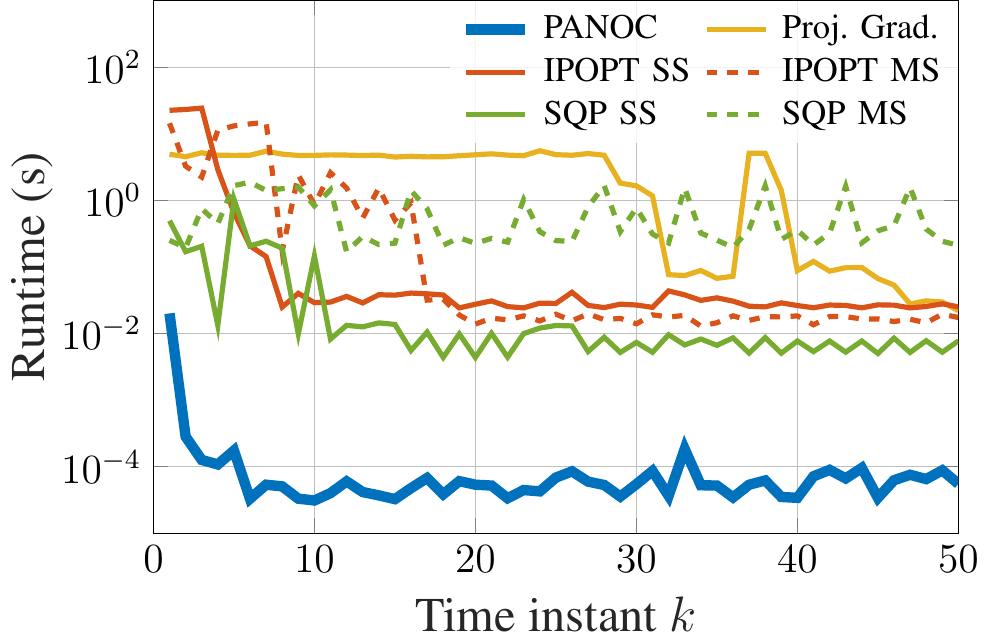}%
 \caption{%
      Comparison of runtime to solve the NMPC problem for four solvers: PANOC, IPOPT and \texttt{fmincon}'s SQP algorithm (for the single and multiple shooting formulations) and projected gradient (ForBES implementation). These timings correspond to the navigation problem presented in Fig.~\ref{fig:obstacle_avoidance_scenarios} (upper right subfigure), starting from the initial point $x_0=(-0.1, -0.2,  \nicefrac{\pi}{5})$ and with $N=50$. The tolerance was set to $3\cdot 10^{-3}$ for all solvers. All solvers are warm-started with their previous solution.
      \label{fig:runtime}
      }%
\end{figure}

The single shooting formulation of Section~\ref{sec:obstacle_avoidance} is solved with PANOC, the interior point solver IPOPT, the forward-backward splitting (FBS) implementation of ForBES and the SQP of MATLAB's \texttt{fmincon}. The problem was also brought in a multiple shooting formulation where obstacle avoidance constraints were imposed as in~\eqref{eq:obst-avoidance-equality} and it was solved using IPOPT and SQP. A comparison of runtime is shown in Fig.~\ref{fig:runtime}. All simulations were executed on an Intel i5-6200U CPU with $\unit[12]{GB}$ RAM machine running Ubuntu 14.04. PANOC exhibits very low runtime and outperforms all other solvers by approximately two orders of magnitude.

\section{Experimental results and discussion}
\begin{figure*}[ht]
\centering
\setlength{\tabcolsep}{0.05em}
\begin{tabular}[t]{rrrrr}
  \includegraphics[width=0.242\textwidth, page=1]{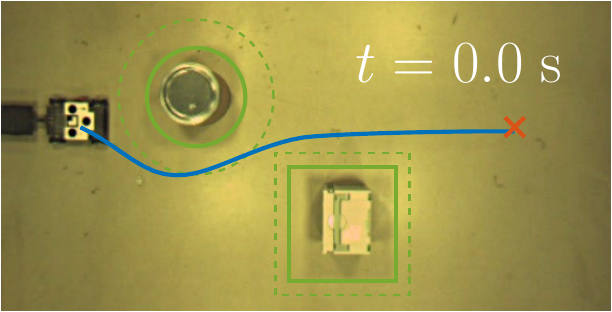} &
  \includegraphics[width=0.242\textwidth, page=2]{snapshots.pdf} &
  \includegraphics[width=0.242\textwidth, page=3]{snapshots.pdf} &
  \includegraphics[width=0.242\textwidth, page=4]{snapshots.pdf}
\end{tabular}
\caption{Point-to-point motion of a holonomic robot with a trailer. The robot avoids a circular and rectangular obstacle indicated with the green lines. The predicted position trajectories are represented in blue and the target position is indicated with the red cross.}
\label{fig:snapshots}
\end{figure*}

The proposed NMPC algorithm is validated on an in-house mobile robot with a trailer (Fig.~\ref{fig:ourbot}). The robot has four independent DC motors with Mecanum wheels, which render it holonomic. A low-level microcontroller implements a velocity controller for each motor. Therefore, from a high-level perspective, the robot can be treated as a velocity-steered holonomic device.
A trailer is attached to the robot and the angle between robot and trailer is measured with a rotary potentiometer. A ceiling camera detects the robot's absolute position and orientation using a marker attached on top of the robot. This information is sent to the robot over Wi-Fi and is merged with local encoder measurements to retrieve a fast and accurate estimate of its pose.
An Odroid XU4 platform runs the NMPC algorithm on board the robotic platform. PANOC is implemented in C following the C89 standard without external dependencies and can, therefore, be readily executed on embedded systems. The C code for performing Alg.~\ref{alg:ad} was generated using the AD tool CasADi \cite{Andersson2013b}.

\begin{figure}
\centering
\includegraphics[width=0.8\linewidth]{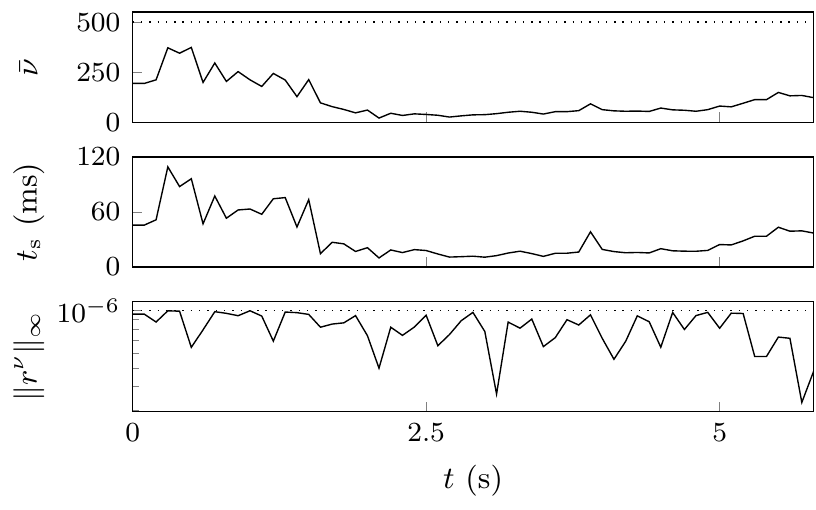}
\caption{Number of iterations $\bar{\nu}$, solve time $t_{s}$ and final residual $\|r^\nu\|_\infty$ for each NMPC cycle.}
\label{fig:iterations}
\end{figure}

The NMPC algorithm is used to steer the vehicle to a desired destination  and trailer orientation $x_{\mathrm{ref}} = (3.77, 1.40, 0.0)$.
The system dynamics is discretized using the Euler method. A quadratic cost function is employed with $Q_k = Q_N = 0.1 \cdot I_3$ and $R_k = 0.01 \cdot I_2$.
Box constraints are imposed on the inputs with $u_{\mathrm{min}} = \unitfrac[-0.8]{m}{s}$ and $u_{\mathrm{max}} = \unitfrac[0.8]{m}{s}$.
The obstacles are modelled as discussed in Section~\ref{sec:reform-obst-avoid} and are slightly enlarged for safety so as to completely avoid collisions. A horizon length of $N{}={}50$ is used and the NMPC control rate is $\unit[10]{Hz}$.
Fig.~\ref{fig:snapshots} shows the resulting motion of the robot and the predicted position trajectories at four time instants. PANOC is terminated when $\|r^\nu\|_{\infty} \leq 10^{-6}$ or if the number of iterations reaches $500$. The number of iterations $\bar{\nu}$, the solving time $t_{\mathrm{s}}$ on the Odroid XU4 and the norm of the fixed point residual at termination at each time instant $k$ are shown in Fig.~\ref{fig:iterations}. 

Additional material and videos from the experiments are found at \href{https://kul-forbes.github.io/PANOC/}{https://kul-forbes.github.io/PANOC}.


\section{Conclusions and research directions}
We proposed a novel framework which enables us to encode nonconvex obstacle avoidance constraints as a smooth nonlinear equality constraint. This offers a  flexible modeling framework and allows the formulation of nonlinear MPC problems with obstacles of complex nonconvex geometry. The resulting MPC problem is solved with PANOC which has favorable theoretical convergence properties and outperforms state-of-the-art NLP solvers. This framework was tested using a C89 implementation of PANOC on a lab-scale  system.

Future work will focus on the development of a proximal Lagrangian framework for the online adaptation of the weight parameters in the obstacle avoidance penalty functions --- cf.~\eqref{eq:soft_penalty} --- so that (predicted) constraint violations, modeled by $\tilde{h}(z_k)$, are below a desired tolerance. Moreover, the use of semismooth Newton directions in PANOC will lead to quadratic convergence and superior performance~\cite{patrinos2014forward}.

\bibliographystyle{IEEEtran}
\bibliography{sources}

\begin{thebibliography}{10}
\providecommand{\url}[1]{#1}
\csname url@samestyle\endcsname
\providecommand{\newblock}{\relax}
\providecommand{\bibinfo}[2]{#2}
\providecommand{\BIBentrySTDinterwordspacing}{\spaceskip=0pt\relax}
\providecommand{\BIBentryALTinterwordstretchfactor}{4}
\providecommand{\BIBentryALTinterwordspacing}{\spaceskip=\fontdimen2\font plus
\BIBentryALTinterwordstretchfactor\fontdimen3\font minus
  \fontdimen4\font\relax}
\providecommand{\BIBforeignlanguage}[2]{{%
\expandafter\ifx\csname l@#1\endcsname\relax
\typeout{** WARNING: IEEEtran.bst: No hyphenation pattern has been}%
\typeout{** loaded for the language `#1'. Using the pattern for}%
\typeout{** the default language instead.}%
\else
\language=\csname l@#1\endcsname
\fi
#2}}
\providecommand{\BIBdecl}{\relax}
\BIBdecl

\bibitem{takahashi1989voronoi}
O.~Takahashi and R.~J. Schilling, ``Motion planning in a plane using
  generalized voronoi diagrams,'' \emph{IEEE Trans. Rob. Autom.}, vol.~5,
  no.~2, pp. 143--150, Apr 1989.

\bibitem{minguez2016motion_planning}
J.~Minguez, F.~Lamiraux, and J.-P. Laumond, \emph{Motion Planning and Obstacle
  Avoidance}.\hskip 1em plus 0.5em minus 0.4em\relax Cham: Springer, 2016, pp.
  1177--1202.

\bibitem{bareiss2015generalized}
D.~Bareiss and J.~van~den Berg, ``Generalized reciprocal collision avoidance,''
  \emph{Int. J. Robot. Res.}, vol.~34, no.~12, pp. 1501--1514, 2015.

\bibitem{wang2014synthesis}
P.~Wang and B.~Ding, ``A synthesis approach of distributed model predictive
  control for homogeneous multi-agent system with collision avoidance,''
  \emph{Int. J. Control}, vol.~87, no.~1, pp. 52--63, 2014.

\bibitem{mercy2017splines}
T.~Mercy, R.~V. Parys, and G.~Pipeleers, ``Spline-based motion planning for
  autonomous guided vehicles in a dynamic environment,'' \emph{IEEE Trans.
  Control Syst. Technol.}, vol.~PP, no.~99, pp. 1--8, 2017.

\bibitem{debrouwere2013time}
F.~Debrouwere, W.~Van~Loock, G.~Pipeleers, M.~Diehl, J.~De~Schutter, and
  J.~Swevers, ``Time-optimal path following for robots with object collision
  avoidance using {Lagrangian} duality,'' in \emph{9th IEEE Int. Workshop Rob.
  Mot. Conntrol (RoMoCo)}, 2013, pp. 186--191.

\bibitem{alrifaee2014collision}
B.~Alrifaee, M.~G. Mamaghani, and D.~Abel, ``Centralized non-convex model
  predictive control for cooperative collision avoidance of networked
  vehicles,'' in \emph{IEEE ISIC}, Oct 2014, pp. 1583--1588.

\bibitem{frasch2013nmpc}
J.~V. Frasch, A.~Gray, M.~Zanon, H.~J. Ferreau, S.~Sager, F.~Borrelli, and
  M.~Diehl, ``An auto-generated nonlinear {MPC} algorithm for real-time
  obstacle avoidance of ground vehicles,'' in \emph{Eur. Control Conf.}, 2013,
  pp. 4136--4141.

\bibitem{turri2013lane}
V.~Turri, A.~Carvalho, H.~E. Tseng, K.~H. Johansson, and F.~Borrelli, ``Linear
  model predictive control for lane keeping and obstacle avoidance on low
  curvature roads,'' in \emph{IEEE Intell. Transp. Syst. Conf.}, Oct 2013, pp.
  378--383.

\bibitem{nocedal2006numerical}
J.~Nocedal and S.~J. Wright, \emph{Numerical Optimization}.\hskip 1em plus
  0.5em minus 0.4em\relax Springer New York, 2006.

\bibitem{boyd2004convex}
S.~Boyd and L.~Vandenberghe, \emph{Convex optimization}.\hskip 1em plus 0.5em
  minus 0.4em\relax Cambridge university press, 2004.

\bibitem{wachter2006intpoint}
A.~W{\"a}chter and L.~T. Biegler, ``On the implementation of an interior-point
  filter line-search algorithm for large-scale nonlinear programming,''
  \emph{Mathem. Progr.}, vol. 106, no.~1, pp. 25--57, Mar 2006.

\bibitem{diehl2005realtime}
M.~Diehl, H.~G. Bock, and J.~P. Schl{\"o}der, ``A real-time iteration scheme
  for nonlinear optimization in optimal feedback control,'' \emph{SIAM J.
  Contr. Optim.}, vol.~43, no.~5, pp. 1714--1736, 2005.

\bibitem{attouch2013fbs}
H.~Attouch, J.~Bolte, and B.~F. Svaiter, ``Convergence of descent methods for
  semi-algebraic and tame problems: proximal algorithms, forward--backward
  splitting, and regularized {Gauss--Seidel} methods,'' \emph{Math. Progr.},
  vol. 137, no.~1, pp. 91--129, Feb 2013.

\bibitem{bertsekas1999nonlinear}
D.~P. Bertsekas, \emph{Nonlinear programming}.\hskip 1em plus 0.5em minus
  0.4em\relax Athena Scientific, 1999.

\bibitem{stella2017simple}
L.~Stella, A.~Themelis, P.~Sopasakis, and P.~Patrinos, ``A simple and efficient
  algorithm for nonlinear model predictive control,'' in \emph{IEEE CDC},
  Melbourne, Australia, 2017.

\bibitem{gilbert1985distance}
E.~Gilbert and D.~Johnson, ``Distance functions and their application to robot
  path planning in the presence of obstacles,'' \emph{IEEE J. Rob. Autom.},
  vol.~1, no.~1, pp. 21--30, Mar 1985.

\bibitem{scheel2000complem}
H.~Scheel and S.~Scholtes, ``Mathematical programs with complementarity
  constraints: Stationarity, optimality, and sensitivity,'' \emph{Math. Op.
  Res.}, vol.~25, no.~1, pp. 1--22, 2000.

\bibitem{liang2012vcc}
Y.-C. Liang and G.-H. Lin, ``Stationarity conditions and their reformulations
  for mathematical programs with vertical complementarity constraints,''
  \emph{J. Optim. Theory Appl.}, vol. 154, no.~1, pp. 54--70, 2012.

\bibitem{dunn1989efficient}
J.~C. Dunn and D.~P. Bertsekas, ``Efficient dynamic programming implementations
  of {N}ewton's method for unconstrained optimal control problems,'' \emph{J.
  Optim. Theory \& Appl.}, vol.~63, no.~1, pp. 23--38, 1989.

\bibitem{patrinos2014forward}
P.~Patrinos, L.~Stella, and A.~Bemporad, ``Forward-backward truncated {N}ewton
  methods for convex composite optimization,'' \emph{arXiv:1402.6655}, Feb.
  2014.

\bibitem{stella2017forward}
L.~Stella, A.~Themelis, and P.~Patrinos, ``Forward-backward quasi-{N}ewton
  methods for nonsmooth optimization problems,'' \emph{Comput. Optim. Appl.},
  vol.~67, no.~3, pp. 443--487, Jul 2017.

\bibitem{themelis2016forward}
A.~Themelis, L.~Stella, and P.~Patrinos, ``Forward-backward envelope for the
  sum of two nonconvex functions: Further properties and nonmonotone
  line-search algorithms,'' \emph{ArXiv preprint}, jun 2016.

\bibitem{PatBem13}
P.~Patrinos and A.~Bemporad, ``Proximal {Newton} methods for convex composite
  optimization,'' in \emph{IEEE CDC}, Dec 2013, pp. 2358--2363.

\bibitem{Andersson2013b}
J.~Andersson, ``A general-purpose software framework for dynamic
  optimization,'' {P}h{D} thesis, KU Leuven, Dept. Electr. Eng. (ESAT/SCD) \&
  Optimization in Engineering Center, October 2013.

\end{thebibliography}

\end{document}